\documentclass[11pt]{article}
\usepackage[ansinew]{inputenc}
\usepackage{array}
\usepackage{color}
\usepackage{amsmath,amsthm,amscd}
\usepackage{amsxtra}
\usepackage{amstext}
\usepackage{amssymb}
\usepackage{hyperref}
\usepackage{latexsym}
\usepackage{amsfonts,cite}
\usepackage{graphics,epstopdf}
\usepackage{epsfig}
\usepackage{amsfonts}
\usepackage{graphicx}
\usepackage{subcaption}
\usepackage{tcolorbox}
\usepackage{array,tabularx,colortbl}
\usepackage[framemethod=TikZ]{mdframed}

\providecommand{\U}[1]{\protect\rule{.1in}{.1in}}

\begin{document}
	
	\sloppy
	\newtheorem{thm}{Theorem}
	\newtheorem{cor}{Corollary}
	\newtheorem{lem}{Lemma}
	\newtheorem{prop}{Proposition}
	\newtheorem{eg}{Example}
	\newtheorem{defn}{Definition}
	\newtheorem{rem}{Remark}
	\newtheorem{note}{Note}
	\numberwithin{equation}{section}
	
	\thispagestyle{empty}
	\parindent=0mm
	
	\begin{center}
		{\large \textbf{A novel advancement in the study of Appell polynomials\\ via Pad\'e rational approximants}}\\ 					
		
		\vspace{0.20cm}
		
		{\bf Giuseppe Dattoli$^{1}$, \bf Subuhi Khan$^{2}$, Ujair Ahmad$^{3}$}\\
		\vspace{0.15cm}
		
		{$^{1}$ENEA - Frascati Research Center, Via Enrico Fermi 45, 00044 Frascati, Rome, Italy.}
		{$^{2,3}$Department of Mathematics, Aligarh Muslim University, Aligarh- 202001, India.}
		\footnote{$^{*}$This work has been done under Junior Research Fellowship (Ref No. 231610072319 dated:29/07/2023) awarded to the third author by University Grand Commission, New Delhi.}
		\footnote{$^{1}$pinodattoli@libero.it (G. Dattoli)}
		\footnote{$^{2}$subuhi2006@gmail.com (Subuhi Khan)}
		\footnote{$^{3}$ujairamu1998@gmail.com (Ujair Ahmad)}
	\end{center}
	\begin{abstract}
		\noindent
	The use of approximants of Padè type are employed to develop 
	a method aimed at opening new perspectives in the theory of 
	Appell polynomials $a_n(x)$, specified by the generating function
	$
	\sum_{n=0}^{\infty} \frac{t^n}{n!} a_n(x) = A(t) e^{xt}.
	$
	In this article, the expansion of amplitude $A(t)$ of the Appell polynomials family in terms of rational approximants yields 
	the possibility of determining the approximation of the $a_n(x)$ in terms of other special polynomials. Application of this approach to Hermite polynomials yields highly accurate approximations in terms of truncated exponential polynomials. Further, monomiality conditions are explored and formalism is extended to consider the Pad\'e approximants within the context of umbral notation.
		
	\end{abstract}
	\parindent=0mm
	\vspace{.25cm}
	
	\noindent
	\textbf{Key Words:}~~Rational approximation; Pad\'e approximation; Appell polynomials; Hermite polynomials; Umbral methods.
	
	\vspace{0.25cm}
	\noindent
	\textbf{2020 Mathematics Subject Classification:}~~ 32E30; 41A21;  11B83;  33C45; 05A40; 41-04.
	
	\section{Introduction}
	The Appell polynomials $A_{n}(x)$, specified by the generating function \cite{APP}
	\begin{equation}\label{paeq1}
		\sum_{n=0}^{\infty}\dfrac{t^{n}}{n!}a_{n}(x) = A(t)e^{xt},
	\end{equation}
	are central focus to our investigation. Here $A(t)$, from now on will be called the amplitude of the Appell polynomials, an infinitely differentiable function defined within a given interval of $t$.\\
	
	This study focuses on developing specific approximations of the amplitude and their implications on the corresponding polynomial sequences. This approach employs Pad\'e approximants \cite{BAK}, which offer a powerful tool for approximating transcendental functions in terms of a rational function \cite{burden}. These approximants realize a generalization of the ordinary Taylor or Maclaurin series expansions.\\
	
	A Pad\'e approximant (PA) of order $m+n$ used to approximate a function $f(x)$ on a closed interval $I$ containing zero is denoted by ${_{[m|n]}f(x)}$. It indicates the following rational function:
	\begin{equation}\label{paeq4}
		{_{[m|n]}f(x)} = \dfrac{P_{m}(x)}{Q_{n}(x)}=\dfrac{\sum_{r=0}^{m}a_{r}x^{r}}{\sum_{r=0}^{n}b_{r}x^{r}}.
	\end{equation}
	For ${_{[m|n]}f(x)}$ to be defined at $x=0$, we require $b_{0} \neq 0$. Without loss of generality, we may assume $b_{0}=1$, since otherwise we can normalize both numerator and denominator by $b_{0}$. The Pad\'e approximation technique chooses $m+n+1$ parameters so that $f^{(k)}(0)={_{[m|n]}f^{(k)}(0)}$, for each $k=0,1,2,\cdots,m+n$.\\
	
	 Suppose $f(x)$ has the Maclaurin series expansion $f(x) := \sum_{r=0}^{\infty}c_{r}x^{r}$. Equating ${_{[m|n]}f(x)}$ to the Maclaurin series expansion of $f(x)$, we have
	\begin{equation}
		\left(\sum_{r=0}^{n}b_{r}x^{r}\right)\left(\sum_{r=0}^{\infty}c_{r}x^{r}\right) \approx \sum_{r=0}^{m}a_{r}x^{r},\quad b_{0}=1,
	\end{equation} 
	so that the rational function for Pad\'e approximation results from the solution of following system of $m+n+1$ linear equations:
	\begin{equation}\label{paeq36}
		\sum_{j=0}^{k}c_{j}b_{k-j} = a_{k}, \quad k = 0,1,2,\cdots,m+n.
	\end{equation}
	The Pad\'e approximants can be utilized to approximate non-rational functions, as for example,  the exponential function $e^{-x}$.\\

	In view of definition \eqref{paeq4}, let us consider the PA ${_{[0|2]}e^{-x}}$:
		\begin{equation}\label{paeq5}
		{_{[0|2]}e^{-x}} = \dfrac{a_{0}}{1+b_{1}x+b_{2}x^{2}}.
		\end{equation}
		Notably, the Maclaurin expansion of order $m$ can be considered as a special case of PA and is denoted as ${_{[m|0]}f(x)}$. Consequently, the second-order expansion of the exponential function $e^{-x}$ writes indeed
		\begin{equation}\label{paeq6}
			e^{-x} \sim 1 - x + \dfrac{x^{2}}{2} = {_{[2|0]}e^{-x}}.
		\end{equation}
	The coefficients are determined by solving the following system of linear equations:
		\begin{align*}
			&b_{2}-b_{1}+\dfrac{1}{2}=0,\\&
		\hspace{0.9cm}	b_{1}-1=0,\\&
		\hspace{1.52cm}	a_{0}=1,
		\end{align*}
		so that, we find
		\begin{equation*}
		a_{0}=b_{1}=1,~~b_{2}=\dfrac{1}{2}.
		\end{equation*}
		Substituting these coefficients in equation \eqref{paeq5}, we find the following second order PA ${_{[0|2]}e^{-x}}$ of exponential function:
		\begin{equation}\label{paeq7}
				{_{[0|2]}e^{-x}} = \dfrac{1}{\left(1 + x + \dfrac{1}{2}x^{2}\right)}.
		\end{equation}
		 Proceeding on the same lines, we obtain the following second order PA ${_{[1|1]}e^{-x}}$ of exponential function:
		\begin{equation}\label{paeq39}
		{_{[1|1]}e^{-x}} = \dfrac{\left(1 - \dfrac{1}{2}x\right)}{\left(1 + \dfrac{1}{2}x\right)}.
		\end{equation}

	Further, we consider the PA ${_{[0|3]}e^{-x}}$:
	\begin{equation}\label{paeq60}
		{_{[0|3]}e^{-x}} = \dfrac{a_{0}}{1 + b_{1}x + b_{2}x^{2} + b_{3}x^{3}}.
	\end{equation}
Since, the third-order Maclaurin expansion of $e^{-x}$ is given by
	\begin{equation}\label{paeq61}
		e^{-x} \sim 1 - x + \dfrac{x^{2}}{2} - \dfrac{x^{3}}{6}.
	\end{equation}
	Therefore Maclaurin expansion \eqref{paeq61} represents PA ${_{[3|0]}e^{-x}}$.\\
	Following the same procedure as above, the PA ${_{[0|3]}e^{-x}}$ is derived as:
	\begin{equation}\label{paeq62}
		{_{[0|3]}e^{-x}} = \dfrac{1}{\left(1 + x + \dfrac{1}{2}x^{2} + \dfrac{1}{6}x^{3}\right)}.
	\end{equation}
In order to demonstrate the comparative accuracy, the comparison between the exponential function with its second order and third order Pad\'e approximants are indicated in figures $1(a)$ and $1(b)$ respectively. Figure 1(b) illustrates that the Pad\'e approximants $[0|3]~\text{and}~[3|0]$ yield better correspondence with the reference function $e^{-x}$ (for further comments, see, \cite{DattoliP}). 
	\begin{figure}[h!]
		\centering
		\begin{subfigure}[b]{0.48\textwidth}
			\includegraphics[width=\textwidth]{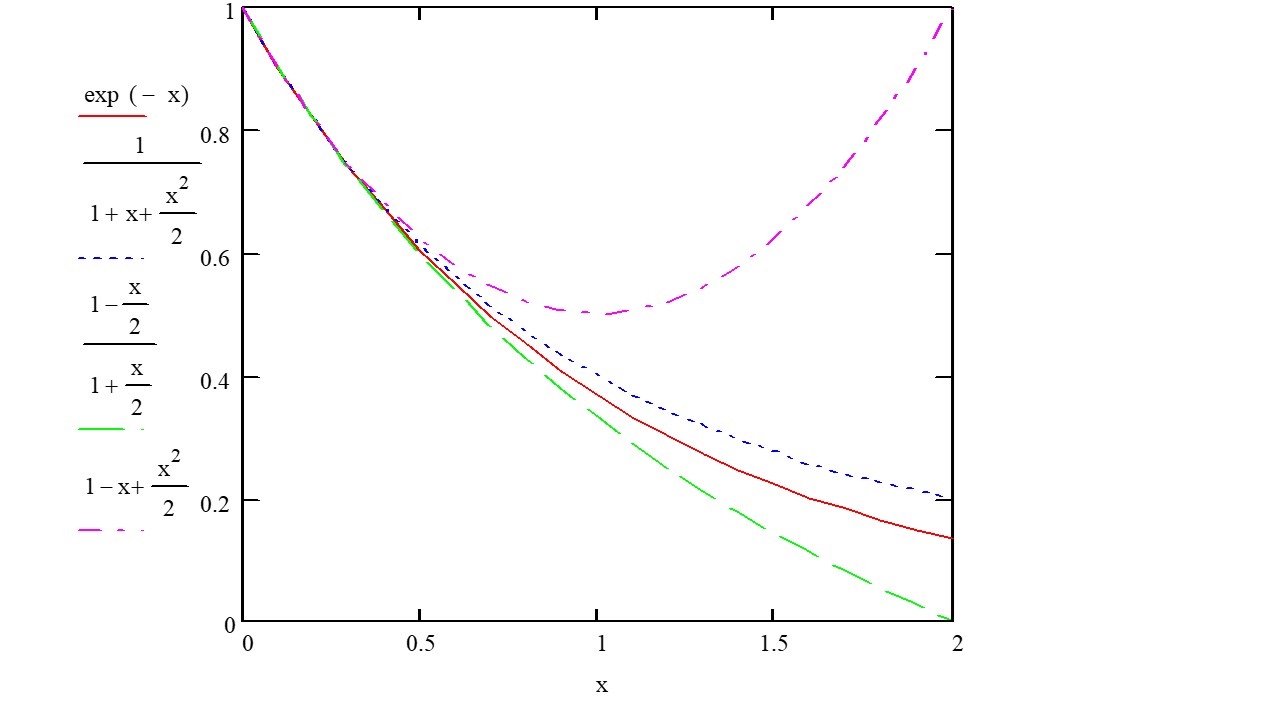}
		\end{subfigure}
		\hfill
		\begin{subfigure}[b]{0.48\textwidth}
			\includegraphics[width=\textwidth]{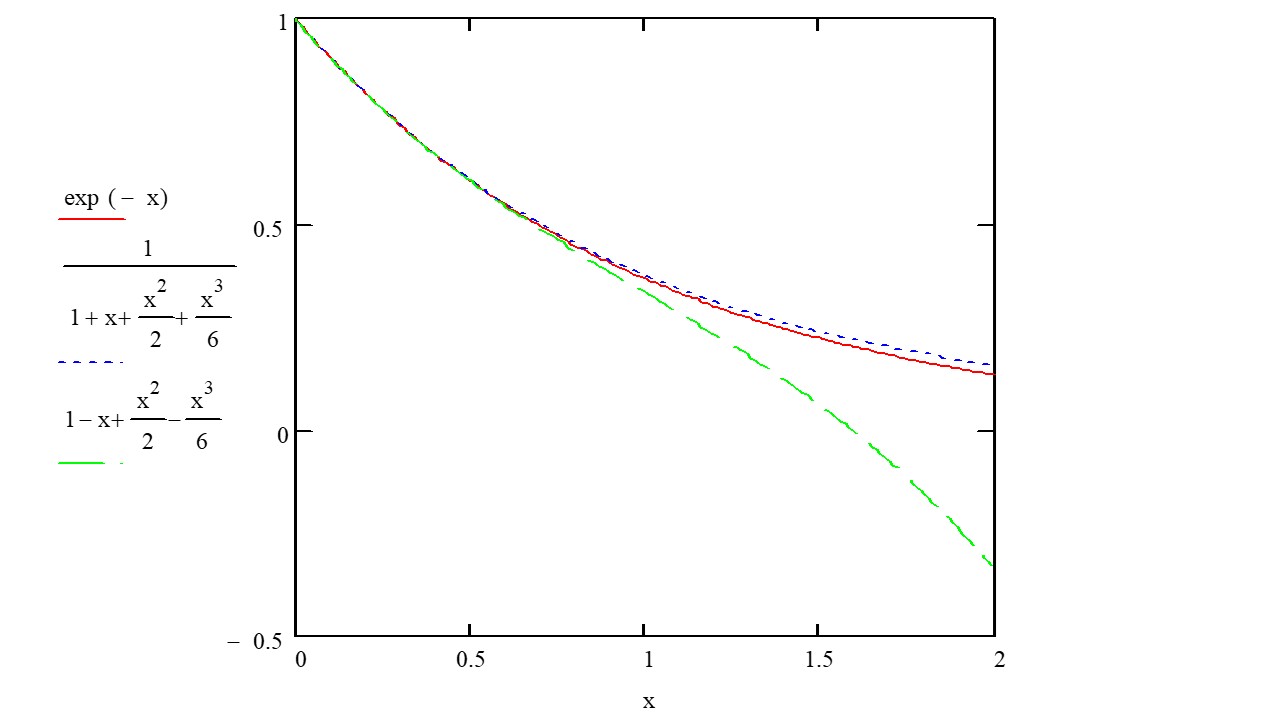}
		\end{subfigure}
		{
			\centering
			\makebox[0.9\textwidth][c]{\em \textbf{Figure 1: Comparison of decaying exponential with Maclaurin}}\\
			\makebox[0.9\textwidth][c]{\em \textbf{and Pad\'e approximations}}\\[0.5em]
			\makebox[0.9\textwidth][c]{\em (a) Red $e^{-x}$, Dot ${_{[0|2]}e^{-x}}$, Slash ${_{[1|1]}e^{-x}}$, Dot Slash ${_{[2|0]}e^{-x}}$}\\
			\makebox[0.9\textwidth][c]{\em (b) Red $e^{-x}$, Dot ${_{[0|3]}e^{-x}}$, Slash ${_{[3|0]}e^{-x}}$}
		}
	\end{figure}\\
	The above example illustrate the process of Padé approximation. Moreover, the accompanying graphs show that Padé approximants provide greater accuracy than the Maclaurin expansions of the same order. In the forthcoming sections, we draw further consequences from the concepts, we have outlined.
\section{Pad\'e approximation of Hermite polynomials}
Hermite polynomials hold an important position across various mathematical domains and can be recognized as a noticeable member of the Appell polynomials family, though this classification requires an appropriate normalization.\\

First, we determine the PA of Hermite polynomials $H_n^{(1)}(x,y)$ of order one defined by the following generating function \cite{BDLS}:
\begin{equation}\label{paeq2}
	\sum_{n=0}^{\infty}\dfrac{t^{n}}{n!}H_n^{(1)}(x,y)=e^{yt} e^{xt}.
\end{equation}
It is to be noted that $H_n^{(1)}(x,y)$ has the following Newton binomials representation:
\begin{equation}\label{paeq3}
	H_n^{(1)}(x,y)=n!\sum_{r=0}^{n}\dfrac{x^{n-r}y^{r}}{(n-r)!r!}=(x+y)^{n}.
\end{equation} 
Since, the PA of amplitude $A(t)$ can be employed to get an approximation of the corresponding special polynomials. Therefore, we have
\begin{equation}\label{paeq8}
	\sum_{n=0}^{\infty}\dfrac{t^{n}}{n!}H_{n}^{(1)}(x,-1)=e^{-t}e^{xt}\sim \sum_{n=0}^{\infty}\dfrac{t^{n}}{n!}\;{_{[r|s]}H_{n}^{(1)}(x,-1)}=\dfrac{P_{r}(t)}{Q_{s}(t)}e^{xt}.
\end{equation}
We note that the polynomials ${_{[r|s]}H_{n}^{(1)}(x,-1)}$ are still belonging to the Appell family and assume that the following relationship holds:
\begin{equation}\label{paeq9}
	H_{n}^{(1)}(x,-1)\sim {_{[r|s]}H_{n}^{(1)}(x,-1)}.
\end{equation}
To verify the accuracy of the previous ansatz, we proceed to determine the second order PA of Hermite polynomials ${_{[1|1]}H_{n}^{(1)}(x,-1)}$ in terms of the first order two-variable truncated exponential polynomials $e_{n}(x,y)$ \cite{DC}.\\

Here, the first order two-variable truncated exponential polynomials $e_{n}(x,y)$ possesses the following series expansion \cite{DC}:
\begin{equation}\label{paeq63}
	e_{n}(x,y)=n!\sum_{r=0}^{n}\dfrac{x^{n-r}y^{r}}{(n-r)!}.
\end{equation}
We know that the amplitude of the truncated exponential polynomials $e_{n}(x,y)$ is
\begin{equation}\label{paeq67}
	A(t)=\dfrac{1}{1-yt}
\end{equation} 
and that they satisfy the following identity:
\begin{equation}\label{paeq64}
	\partial_{x} e_{n}(x,y)=ne_{n-1}(x,y).
\end{equation}
		We recall the following operational representation of Appell polynomials\cite{DG}:
		\begin{equation}\label{paeq10}
			a_{n}(x)=A(\partial_{x})\{x^{n}\}.
		\end{equation}
			Since, $A(t)=e^{-t}$ for the Hermite polynomials $H_{n}^{(1)}(x,-1)$, therefore using approximation expression \eqref{paeq39} in equation \eqref{paeq10}, we obtain the following second order Pad\'e approximant of Hermite polynomials ${_{[1|1]}H_{n}^{(1)}(x,-1)}$:
			\begin{equation}\label{paeq11}
				{_{[1|1]}H_{n}^{(1)}(x,-1)}={_{[1|1]}e^{-\partial _{x}}}\{x^{n}\}=\dfrac{\left(1-\dfrac{1}{2} \partial_{x}\right)}{\left(1+\dfrac{1}{2} \partial_{x}\right)}\{x^{n}\},
			\end{equation}
		Further, we derive the explicit representation of the second order PA of Hermite polynomials ${_{[1|1]}H_{n}^{(1)}(x,-1)}$. Rewriting result \eqref{paeq11} in the following form:
		\begin{equation}\label{paeq12}
			{_{[1|1]}H_{n}^{(1)}(x,-1)}=	\dfrac{\left(1-\dfrac{1}{2} \partial_{x}\right)}{\left(1+\dfrac{1}{2} \partial_{x}\right)}\{x^{n}\}=\left(\left(1-\dfrac{1}{2}\partial_{x}\right)\sum_{r=0}^{\infty}\left(-\dfrac{1}{2}\partial_{x}\right)^{r}\right)\{x^{n}\}.
		\end{equation}
In view of equations \eqref{paeq12}, \eqref{paeq63} and \eqref{paeq64}, we find the following explicit representation of PA ${_{[1,1]}H_{n}^{(1)}(x,-1)}$ in terms of $e_{n}(x,y)$:
	\begin{equation}\label{paeq13}
	{_{[1|1]}H_{n}^{(1)}(x,-1)}=e_{n}\left(x,-\dfrac{1}{2}\right)-\dfrac{n}{2}e_{n-1}\left(x,-\dfrac{1}{2}\right),
\end{equation}
		Since, $A(t)=e^{-t}$ is the amplitude of $H_{n}^{(1)}(x,-1)$, therefore, by using the same procedure, it follows that 
		\begin{equation*}
			{_{[2|1]}e^{-t}}=\dfrac{\left(1-\dfrac{2}{3}t+\dfrac{1}{6}t^{2}\right)}{\left(1+\dfrac{1}{3}t\right)}.
		\end{equation*}
		Now, using operational definition \eqref{paeq10}, we find the third order PA ${_{[2|1]}H_{n}^{(1)}(x,-1)}$ of Hermite polynomials as:
		\begin{equation}\label{paeq14}
			{_{[2|1]}H_{n}^{(1)}(x,-1)}=\dfrac{\left(1-\dfrac{2}{3}\partial_{x}+\dfrac{1}{6}\partial_{x}^{2}\right)}{\left(1+\dfrac{1}{3}\partial_{x}\right)}\{x^{n}\}.
		\end{equation}
	Following the same procedure, sketched before, we eventually obtain the following explicit representation of third order PA ${_{[2|1]}H_{n}^{(1)}(x,-1)}$ in terms of $e_{n}(x,y)$:
		\begin{equation}\label{paeq15}
			{_{[2|1]}H_{n}^{(1)}(x,-1)}=\dfrac{P_{2}(\partial_{x})}{Q_{1}(\partial_{x})}\{x^{n}\}=e_{n}\left(x,-\dfrac{1}{3}\right)-\dfrac{2}{3}n\;e_{n-1}\left(x,-\dfrac{1}{3}\right)+\dfrac{1}{6}n(n-1)e_{n-2}\left(x,-\dfrac{1}{3}\right).
		\end{equation}
	
	The same procedure can be adopted to obtain the higher order PA ${_{[r|s]}H_{n}^{(1)}(x,-1)}$ for $H_{n}^{(1)}(x,-1)$.\\
	
In Figures $2(a)-2(c)$, we have shown a comparison between $H_n^{(1)} (x,-1)$ with its PA ${_{[1|1]}H_{n}^{(1)}(x,-1)}$ for specific values of index. It indicates that the approximation becomes less accurate with increasing $n$. In Figure $2(d)$ a net improvement is shown with the up grading of the approximant to ${_{[2|1]}H_{n}^{(1)}(x,-1)}$.\\
	
	\begin{figure}[h!]
		\centering
		\begin{subfigure}[b]{0.48\textwidth}
			\includegraphics[width=\textwidth]{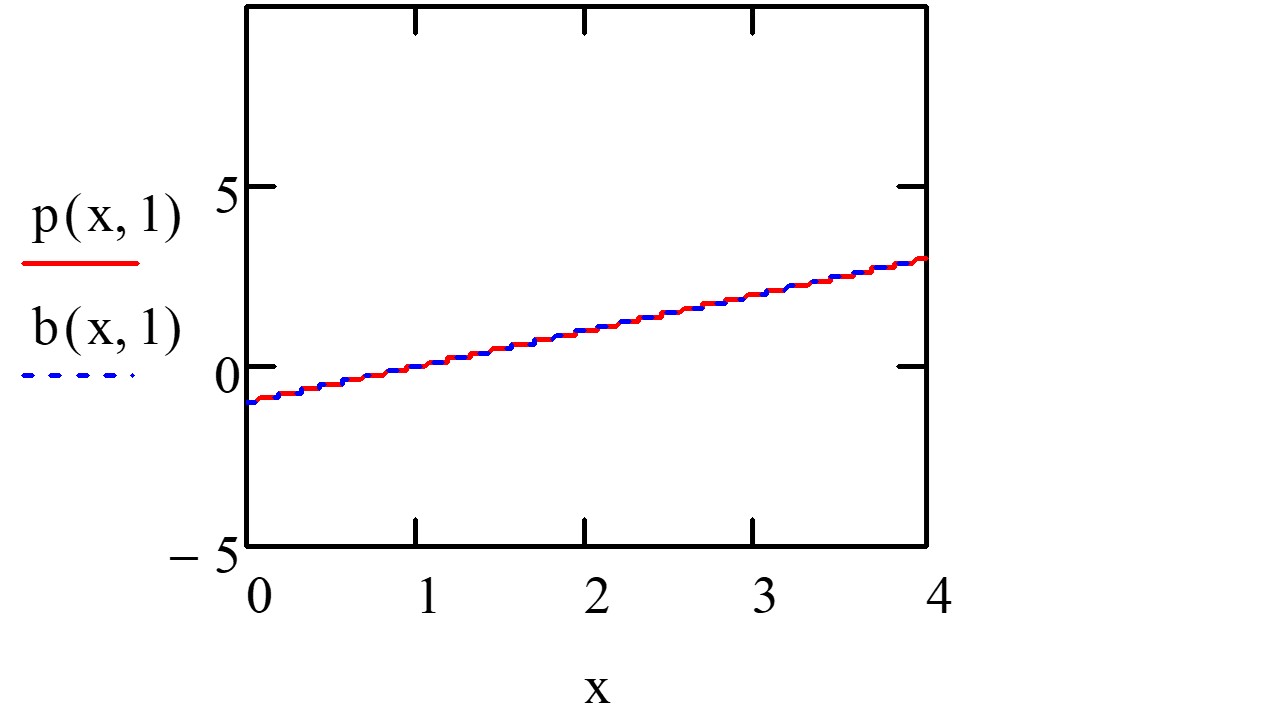}
			\caption*{\centering{\emph{Figure 2(a): Red $p(x,1)= H_{1}^{(1)}(x,-1)$, Dot $b(x,1)={_{[1|1]}H_{1}^{(1)}(x,-1)}$}}}
		\end{subfigure}
		\hfill
		\begin{subfigure}[b]{0.48\textwidth}
			\includegraphics[width=\textwidth]{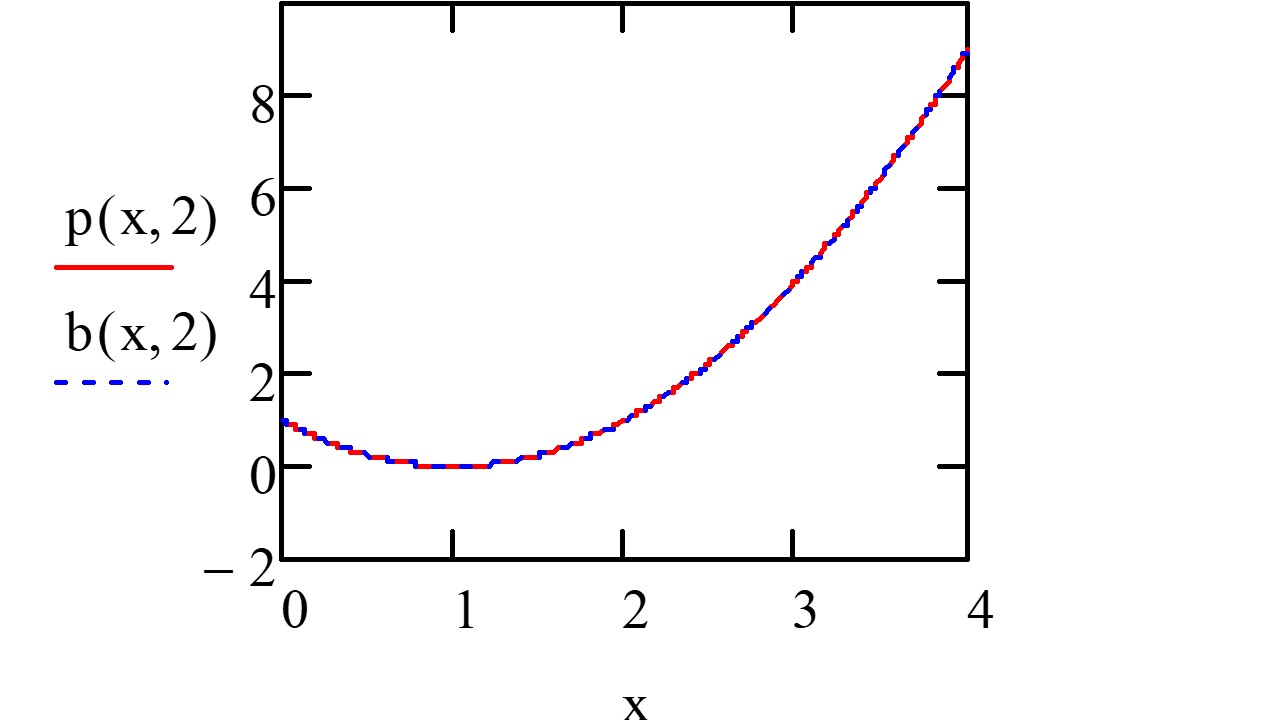}
			\caption*{\centering{\emph{figure 2(b): Red $p(x,2)=H_{2}^{(1)}(x,-1)$, Dot $b(x,2)={_{[1|1]}H_{2}^{(1)}(x,-1)}$}}}
		\end{subfigure}
	\end{figure}
	\begin{figure}[h!]
		\centering
		\begin{subfigure}[b]{0.48\textwidth}
			\includegraphics[width=\textwidth]{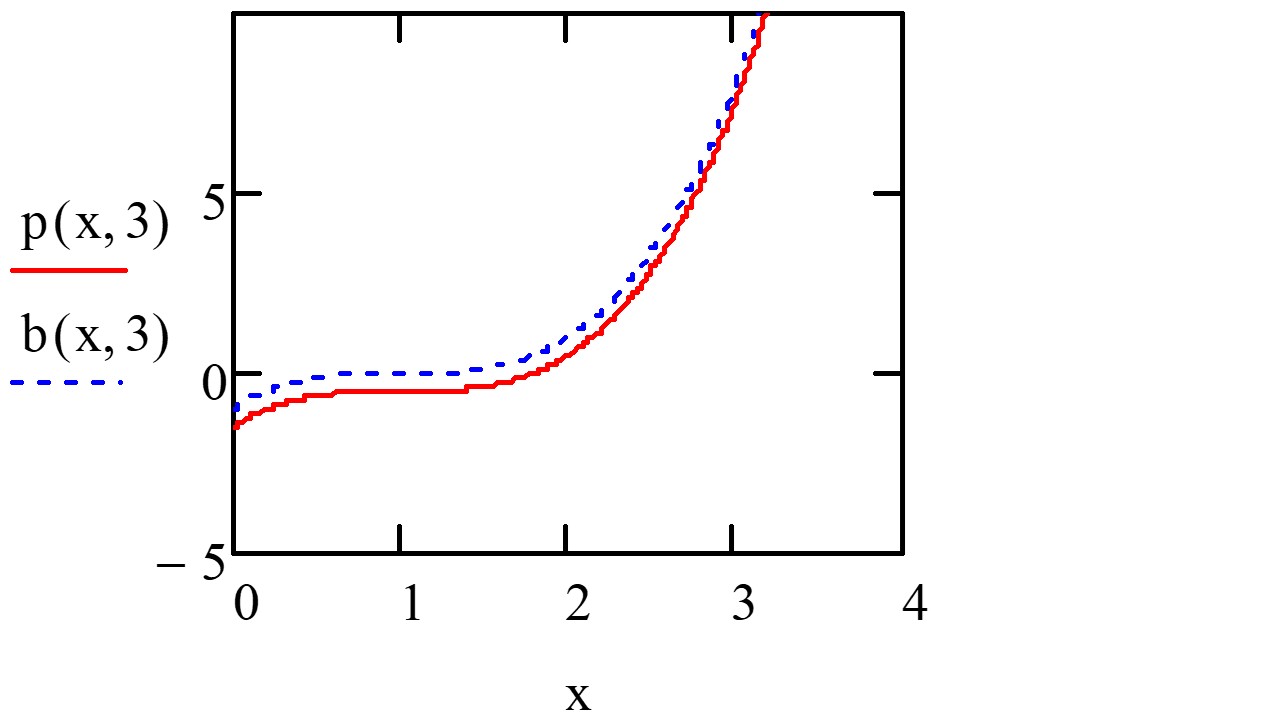}
			\caption*{\centering{\emph{Figure 2(c): Red $p(x,3)=H_{3}^{(1)}(x,-1)$, Dot $b(x,3)={_{[1|1]}H_{3}^{(1)}(x,-1)}$}}}
		\end{subfigure}
		\hfill
		\begin{subfigure}[b]{0.48\textwidth}
			\includegraphics[width=\textwidth]{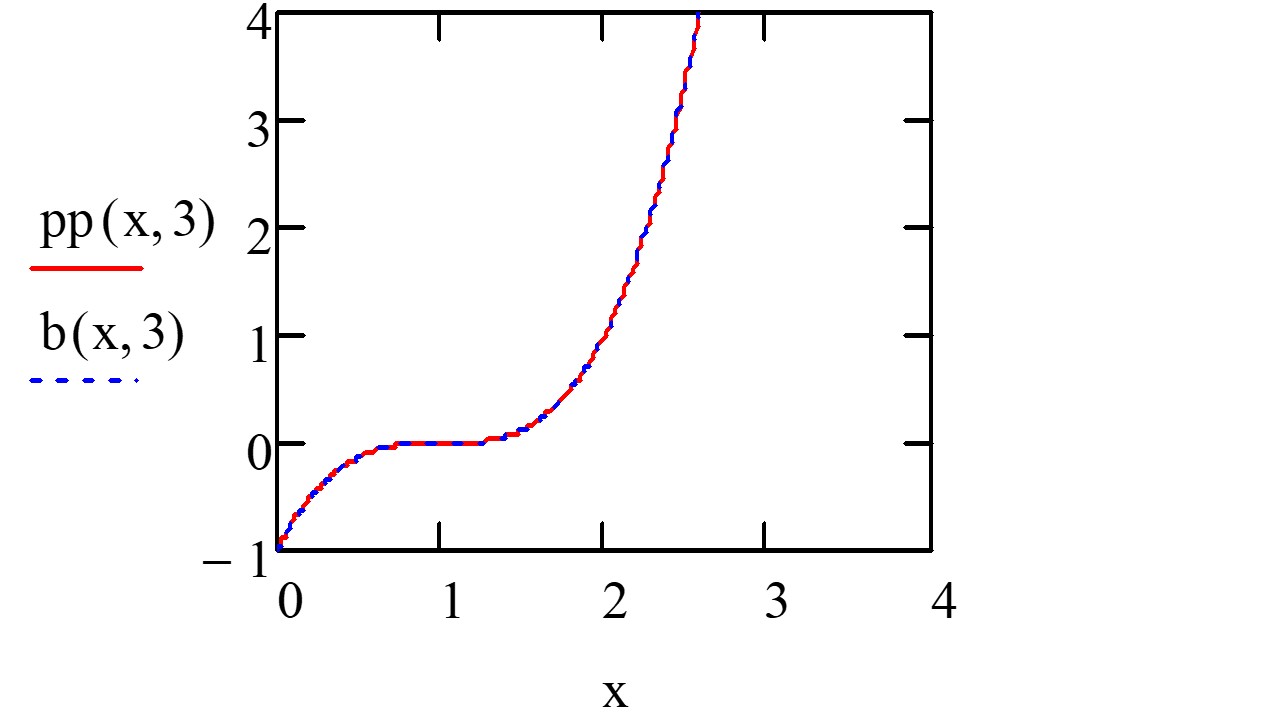}
			\caption*{\centering{\emph{figure 2(d): Red $pp(x,3)=H_{3}^{(1)}(x,-1)$, Dot $b(x,3)={_{[2|1]}H_{3}^{(1)}(x,-1)}$}}}
		\end{subfigure}
	\end{figure}
	Now, we extend the concepts, we have outlined to the two-variable Hermite polynomials \cite{Appell} defined by the following series:
	\begin{equation}\label{paeq22}
		H_{n}^{(2)}(x,y)=n!\sum_{r=0}^{\small\left\lceil\frac{n}{2}\right\rceil}\dfrac{x^{n-2r}y^{r}}{(n-2r)!r!},
	\end{equation}
	which are specified by the Appell amplitude $e^{yt^{2}}$ and the relevant PA are obtained along the same lines adopted for the binomial polynomials \eqref{paeq3}. The ordinary Hermite polynomials $He_{n}(x)$ are a particular case of the two-variable Hermite polynomials and are recognized as $H_{n}^{(2)}\left(x,-\dfrac{1}{2}\right)$.\\
	
	 Next, we establish an identity for the second order PA ${_{[1|1]}He_{n}(x)}$ of ordinary Hermite polynomials in terms of the second order two-variable truncated polynomials $e_{n}^{(2)}(x,y)$ defined by
	\begin{equation}\label{paeq55}
		e_{n}^{(2)}(x,y)=n!\sum_{r=0}^{\small\left\lceil\frac{n}{2}\right\rceil}\dfrac{x^{n-2r}y^{r}}{(n-2r)!}
	\end{equation}
	in the form of following result:
	\begin{thm}
		For the second order Pad\'e approximant of ordinary Hermite polynomials ${_{[1|1]}He_{n}(x)}$, the following explicit representation in terms of second order two-variable truncated exponential polynomials $e_{n}^{(2)}(x,y)$ holds true:
		\begin{equation}\label{paeq23}
			{_{[1|1]}He_{n}(x)}=\dfrac{\left(1-\dfrac{1}{4} \partial_{x}^{2}\right)}{\left(1+\dfrac{1}{4} \partial_{x}^{2}\right)}\{x^{n}\}=e_{n}^{(2)}\left(x,-\dfrac{1}{4}\right)-\dfrac{n(n-1)}{4}e_{n-2}^{(2)}\left(x,-\dfrac{1}{4}\right).
		\end{equation}
		\begin{proof}
			Replacing $x$ by $t^{2}/2$ in approximation equation \eqref{paeq39}, we obtain the second order PA ${_{[1|1]}}A(t)$ of amplitude $A(t)=e^{-t^{2}/2}$:
			\begin{equation*}
				{_{[1|1]}}e^{-t^{2}/2}=\dfrac{\left(1-\dfrac{1}{4} t^{2}\right)}{\left(1+\dfrac{1}{4} t^{2}\right)}.
			\end{equation*}
			In view of operational definition \eqref{paeq10}, we have
			\begin{equation*}
				{_{[1|1]}He_{n}(x)}={_{[1|1]}e^{-\partial_{x}^{2}/2}}{x^{n}}=\dfrac{\left(1-\dfrac{1}{4} \partial_{x}^{2}\right)}{\left(1+\dfrac{1}{4} \partial_{x}^{2}\right)}\{x^{n}\}.
			\end{equation*}
			Rewrite the above equation as:
			\begin{equation*}
				{_{[1|1]}He_{n}(x)}=\left(\left(1-\dfrac{1}{4}\partial_{x}^{2}\right)\sum_{r=0}^{\infty}\left(-\dfrac{1}{4}\partial_{x}^{2}\right)^{r}\right)\{x^{n}\},
			\end{equation*}
			which after simplification yields desired assertion \eqref{paeq23}.
		\end{proof}
	\end{thm}
	
The approximation \eqref{paeq23} is fairly good and is significantly improved by the inclusion of further terms in the Pad\'e expansion (see, Figure $3$).\\
	\begin{figure}[h!]
		\centering
		\begin{subfigure}[b]{0.48\textwidth}
			\includegraphics[width=\textwidth]{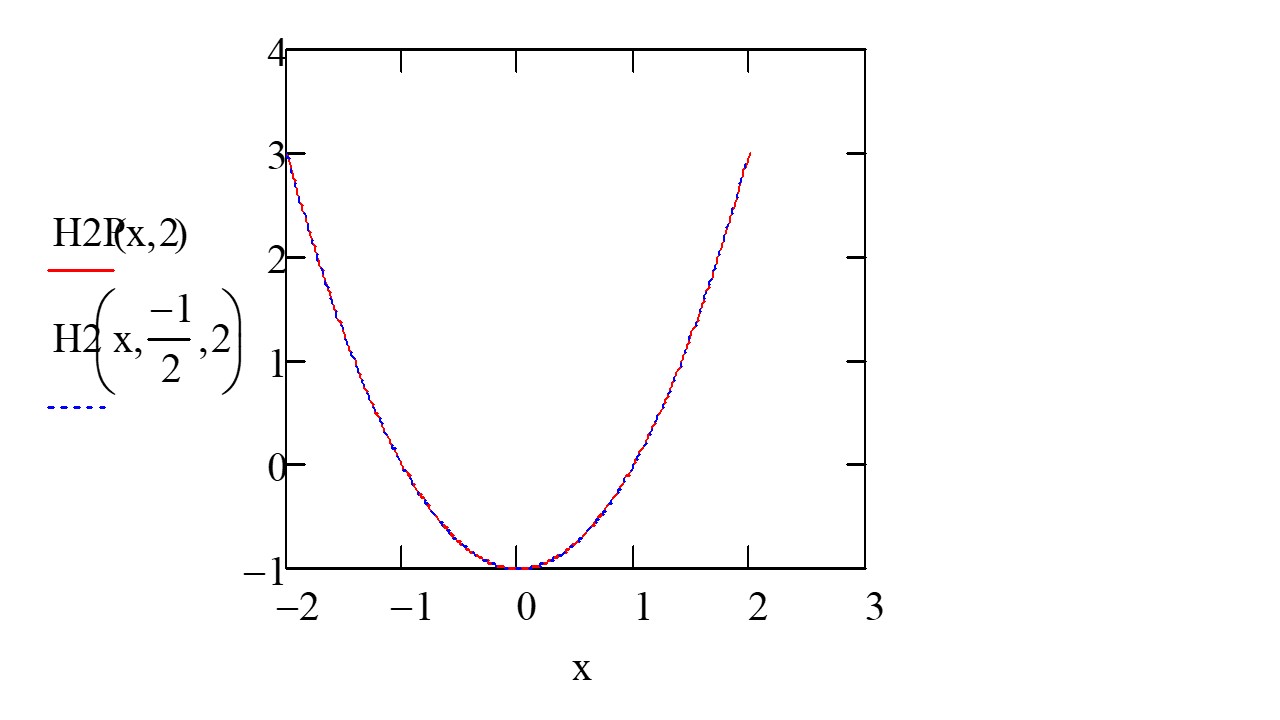}
			\caption*{\centering{\emph{Figure 3(a): Red $H2P(x,2)= {_{[1|1]}He_{2}(x)}$, Dot $H2\left(x,\dfrac{-1}{2},2\right)=He_{2}(x)$}}}
		\end{subfigure}
		\hfill
		\begin{subfigure}[b]{0.48\textwidth}
			\includegraphics[width=\textwidth]{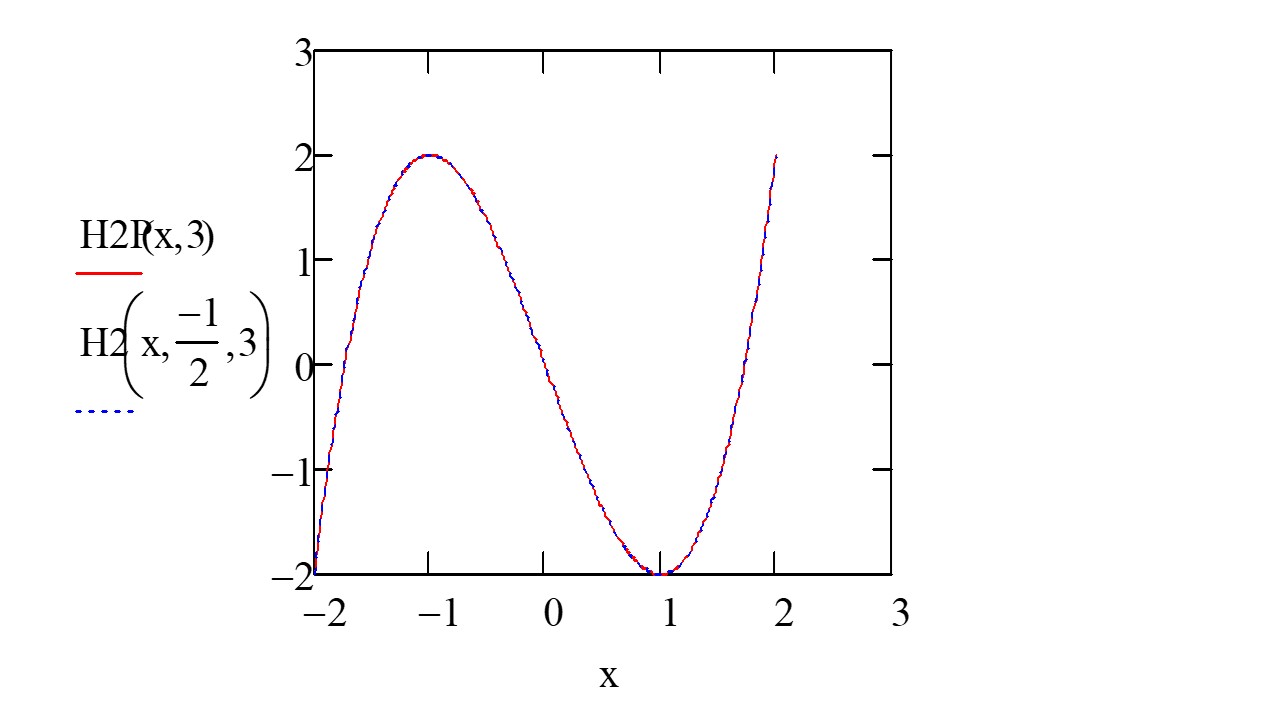}
			\caption*{\centering{\emph{figure 3(b): Red $H2P(x,3)= {_{[1|1]}He_{3}(x)}$, Dot $H2\left(x,\dfrac{-1}{2},3\right)=He_{3}(x)$}}}
		\end{subfigure}
	\end{figure}
	\begin{figure}[h!]
		\centering
		\begin{subfigure}[b]{0.48\textwidth}
			\includegraphics[width=\textwidth]{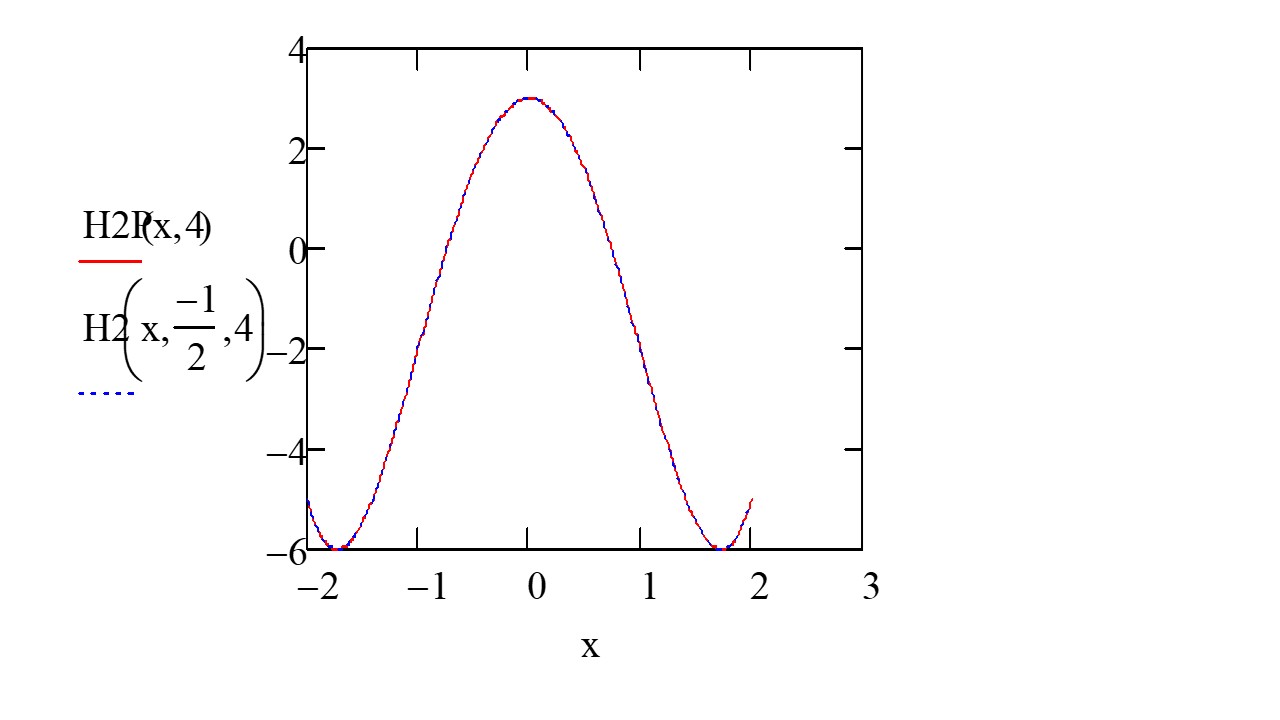}
			\caption*{\centering{\emph{Figure 3(c): Red $H2P(x,4)= {_{[1|1]}He_{4}(x)}$, Dot $H2\left(x,\dfrac{-1}{2},4\right)=He_{4}(x)$}}}
		\end{subfigure}
		\hfill
		\begin{subfigure}[b]{0.48\textwidth}
			\includegraphics[width=\textwidth]{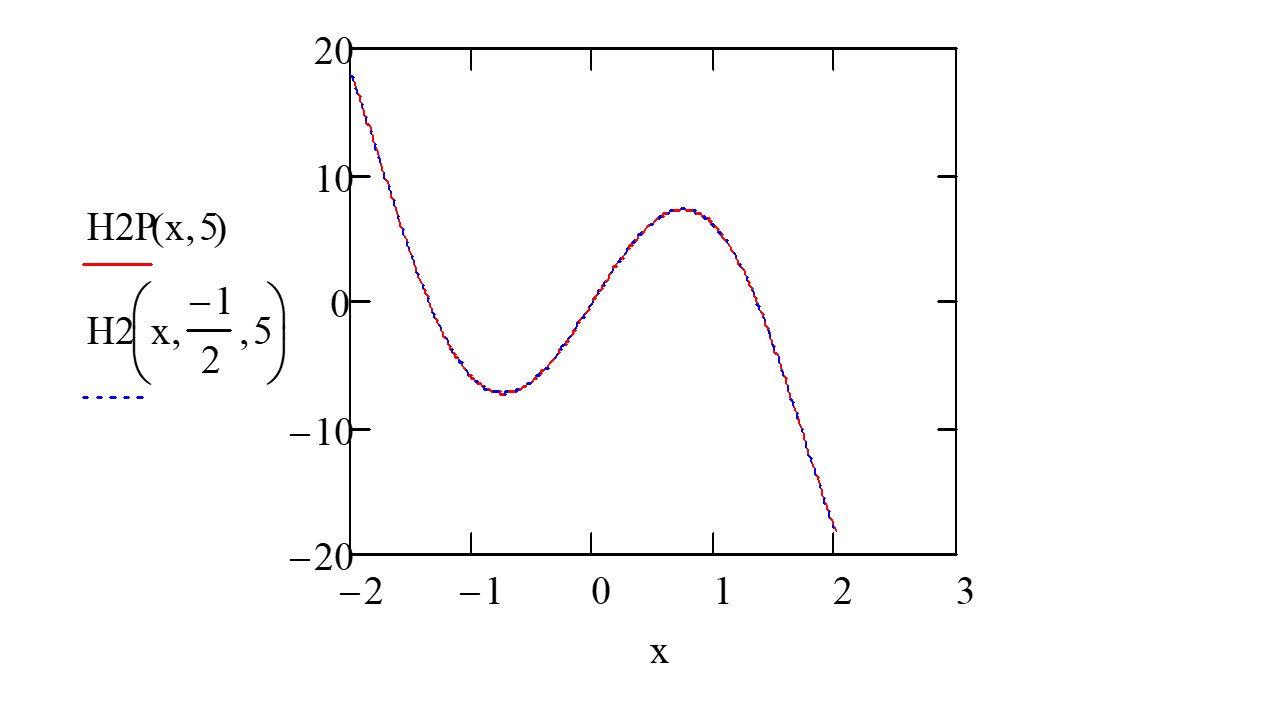}
			\caption*{\centering{\emph{figure 3(d): Red $H2P(x,5)= {_{[1|1]}He_{5}(x)}$, Dot $H2\left(x,\dfrac{-1}{2},5\right)=He_{5}(x)$}}}
		\end{subfigure}
	\end{figure}
	
	\begin{figure}[h!]
		\centering
		\begin{subfigure}[b]{0.48\textwidth}
			\includegraphics[width=\textwidth]{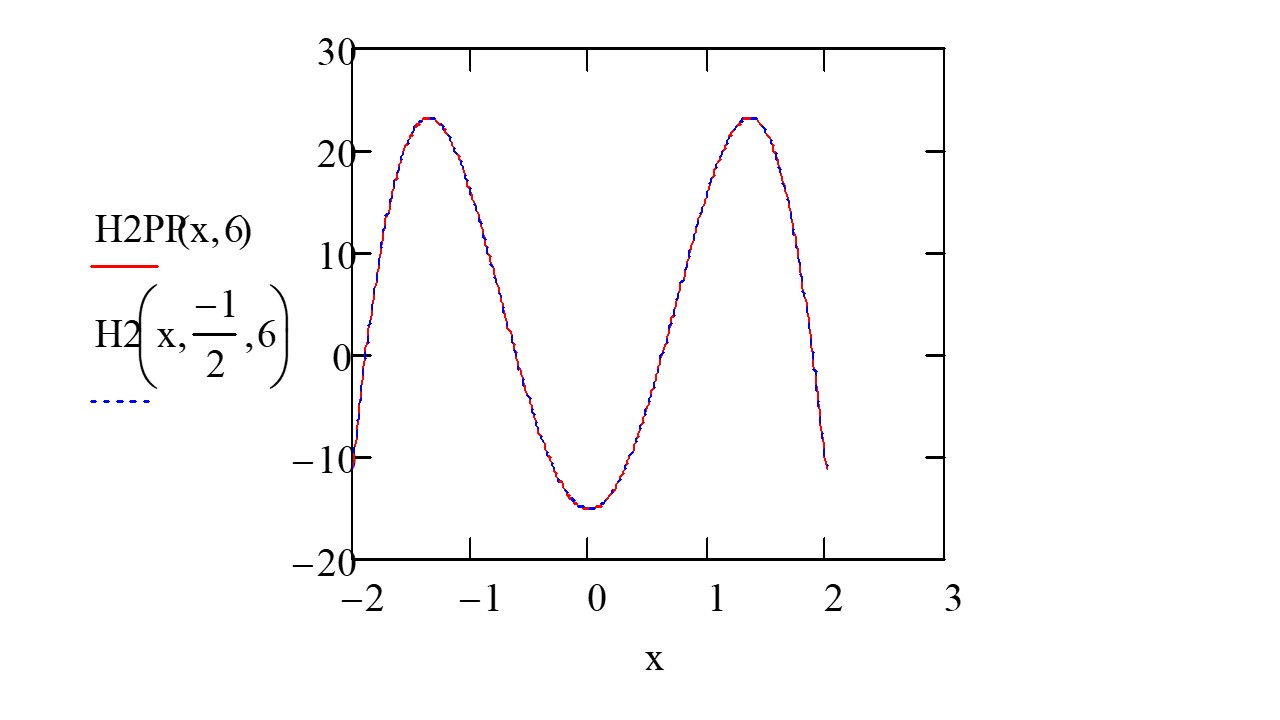}
			\caption*{\centering{\emph{Figure 3(e): Red $H2P(x,6)= {_{[1|1]}He_{6}(x)}$, Dot $H2\left(x,\dfrac{-1}{2},6\right)=He_{6}(x)$}}}
		\end{subfigure}
		\hfill
		\begin{subfigure}[b]{0.48\textwidth}
			\includegraphics[width=\textwidth]{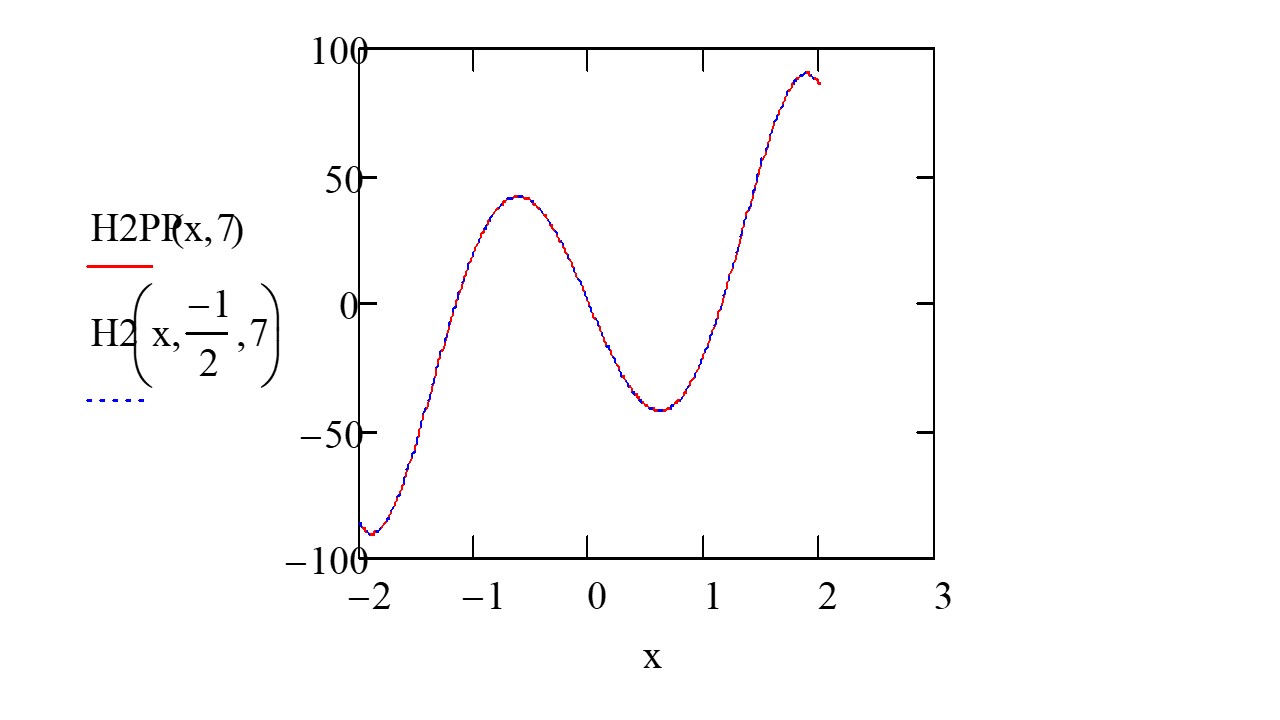}
			\caption*{\centering{\emph{figure 3(f): Red $H2P(x,7)= {_{[1|1]}He_{7}(x)}$, Dot $H2\left(x,\dfrac{-1}{2},7\right)=He_{7}(x)$}}}
		\end{subfigure}
		\centering
		\makebox[0.9\textwidth][c]{\em \textbf{Figure 3: Hermite polynomials $He_{n}(x)$ (dot) and their}}\\
		\makebox[0.9\textwidth][c]{\em \textbf{Pad\'e Approximation (red) for $n=2,3,\dots,7$.}}
	\end{figure}

	Before closing this section, we consider some interesting results in terms of the two-variable Chebyshev polynomials of second kind $U_{n}(x,y)$ \cite{SLicciardi}, defined by the following series expansion and generating function:
	\begin{equation*}
		U_{n}(x,y)=(-1)^{n}\sum_{r=0}^{\small\left\lceil\frac{n}{2}\right\rceil}\dfrac{(n-r)!(x)^{n-2r}(-y)^{r}}{(n-2r)!r!}
	\end{equation*}
	and
	\begin{equation}\label{paeq56}
		\dfrac{1}{\left(1+xt+yt^{2}\right)}=\sum_{n=0}^{\infty}U_{n}(x,y)t^{n},
	\end{equation}
	respectively.\\ Here, by taking $-2x$ at the place of $x$ and $y=1$, the ordinary second kind Chebyshev polynomials $U_{n}(x)$ \cite{AND} are obtained.\\
	
	Further, taking $x=a,~y=b$ in generating function \eqref{paeq56}, we have
		\begin{equation}\label{paeq57}
			\dfrac{1}{\left(1+at+bt^{2}\right)}=\sum_{n=0}^{\infty}U_{n}(a,b)t^{n}.
		\end{equation}

	An interesting example is provided by using the approximant ${_{[0|2]}e^{-t}}$ in equation \eqref{paeq7}, which in the case of exponential function reads
	\begin{equation}\label{paeq24}
		{_{[0|2]}e^{-t}}=\dfrac{1}{\left(1+t+\dfrac{1}{2}t^{2}\right)}=\sum_{n=0}^{\infty}U_{n}\left(1,\dfrac{1}{2}\right)t^{n}, \quad |t|<1.
	\end{equation}
	Accordingly, the second order PA ${_{[0|2]}He_{n}(x)}$ of Hermite polynomials is expressed in terms of the Chebyshev polynomials $U_{n}(x,y)$ in the following result:
	\begin{thm}
		For the second order Pad\'e approximant ${_{[0|2]}He_{n}(x)}$ of ordinary Hermite polynomials, the following expression in terms of two-variable Chebyshev polynomials of second kind $U_{n}(x,y)$ holds true:
		\begin{equation}\label{paeq25}
			{_{[0|2]}He_{n}(x)}=n!\sum_{r=0}^{\small\left\lceil\frac{n}{2}\right\rceil}\dfrac{U_{r}\left(\dfrac{1}{2},\dfrac{1}{8}\right)x^{n-2r}}{(n-2r)!}.
		\end{equation}
		\begin{proof}
		Replacing $x$ by $t^{2}/2$ in approximation equation \eqref{paeq7}, we have the following second order PA ${_{[0|2]}A(t)}$ of amplitude $A(t)=e^{-t^{2}/2}$:
		\begin{equation*}
			{_{[0|2]}}e^{-t^{2}/2}=\dfrac{1}{\left(1+\dfrac{1}{2} t^{2}+\dfrac{1}{8}t^{4}\right)}.
		\end{equation*}
	Making use of operational definition \eqref{paeq10} in the above equation, it follows that
			\begin{equation}\label{paeq58}
				{_{[0|2]}He_{n}(x)}=\dfrac{1}{\left(1+\dfrac{1}{2}\partial_{x}^{2}+\dfrac{1}{8}\partial_{x}^{4}\right)}\{x^{n}\}.
			\end{equation}
			Now, in view of equation \eqref{paeq57}, we have
			\begin{equation*}
				\dfrac{1}{\left(1+\dfrac{1}{2}t+\dfrac{1}{8}t^{2}\right)}=\sum_{r=0}^{\infty}U_{r}\left(\dfrac{1}{2},\dfrac{1}{8}\right)t^{r},
			\end{equation*}
		which on multiplying both side by $x^{n}$ and replacing $t$ by $\partial_{x}^{2}$ takes the form
			\begin{equation}\label{paeq73}
				\dfrac{1}{\left(1+\dfrac{1}{2}\partial_{x}^{2}+\dfrac{1}{8}\partial_{x}^{4}\right)}\{x^{n}\}=\sum_{r=0}^{\infty}U_{r}\left(\dfrac{1}{2},\dfrac{1}{8}\right)\partial_{x}^{2r}\{x^{n}\}.
			\end{equation}
			Comparing the above equation with equation \eqref{paeq58} and simplifying, assertion \eqref{paeq25} is obtained.
		\end{proof}
	\end{thm}
It is important to observe that approximation \eqref{paeq25} in the range $|x|<1$ for the first five $He_{n}(x)$ is extremely good and better approximations for larger $n$-values. The discrepancy is easily recovered by adding further terms in the corresponding Pad\'e expansion.\\
	
	The higher order approximants can be obtained along the same lines. The representation of fifth order PA ${_{[3|2]}He_{n}(x)}$ of ordinary Hermite polynomials in terms of Chebyshev polynomials $U_{n}(x,y)$ is established in the next result:
	\begin{thm}
		 For the fifth order Pad\'e approximant ${_{[3|2]}He_{n}(x)}$ of ordinary Hermite polynomials, the following representation holds true:
		\begin{equation}\label{paeq28}
			{_{[3|2]}He_{n}(x)}=d_{n}(x)-\dfrac{3}{10}\dfrac{n!}{(n-2)!}d_{n-2}(x)+\dfrac{3}{80}\dfrac{n!}{(n-4)!}d_{n-4}(x)-\dfrac{1}{480}\dfrac{n!}{(n-6)!}d_{n-6}(x),
		\end{equation}
		where
		\begin{equation}\label{paeq27}
		d_{n}(x)	=n!\sum_{r=0}^{\small\left\lceil\frac{n}{2}\right\rceil}\dfrac{U_{r}\left(\dfrac{1}{5},\dfrac{1}{80}\right)x^{n-2r}}{(n-2r)!}.
		\end{equation}
		\begin{proof}
			Using the method outline in introduction, the fifth order PA ${_{[3|2]}e^{-t^{2}/2}}$ of amplitude $A(t)=e^{-t^{2}/2}$ is expressed as:
			\begin{equation}\label{paeq79}
				{_{[3|2]}}A(t)={_{[3|2]}}e^{-t^{2}/2}=\dfrac{\left(1-\dfrac{3}{10}t^{2}+\dfrac{3}{80}t^{4}-\dfrac{1}{480}t^{6}\right)}{\left(1+\dfrac{1}{5}t^{2}+\dfrac{1}{80}t^{4}\right)}.
			\end{equation}
			Making use of approximation equation \eqref{paeq79} in operational definition \eqref{paeq10}, we obtain the following fifth order PA ${_{[3|2]}He_{n}(x)}$ of ordinary Hermite polynomials:
			\begin{equation}\label{paeq26}
				{_{[3|2]}He_{n}(x)}=\dfrac{\left(1-\dfrac{3}{10}\partial_{x}^{2}+\dfrac{3}{80}\partial_{x}^{4}-\dfrac{1}{480}\partial_{x}^{6}\right)}{\left(1+\dfrac{1}{5}\partial_{x}^{2}+\dfrac{1}{80}\partial_{x}^{4}\right)}\{x^{n}\}=\dfrac{P_{3}(\partial_{x}^{2})}{Q_{2}(\partial_{x}^{2})}\{x^{n}\}.
			\end{equation}
			In view of equation \eqref{paeq73}, we have
			\begin{equation}\label{paeq74}
				\dfrac{1}{Q_{2}(\partial_{x}^{2})}\{x^{n}\}=\dfrac{1}{\left(1+\dfrac{1}{5}\partial_{x}^{2}+\dfrac{1}{80}\partial_{x}^{4}\right)}\{x^{n}\}=n!\sum_{r=0}^{\small\left\lceil\frac{n}{2}\right\rceil}\dfrac{U_{r}\left(\dfrac{1}{5},\dfrac{1}{80}\right)x^{n-2r}}{(n-2r)!}.
			\end{equation}
			Using notation \eqref{paeq27} in equation \eqref{paeq74} and taking the actions of $P_{3}(\partial_{x}^{2})$ on the resultant equation, we get assertion \eqref{paeq28}.
		\end{proof}
	\end{thm}
	The comparison between $He_{n}(x)$ and its approximated forms is given in Figure $4$ where we have reported \eqref{paeq28}.\\
	\begin{figure}[h!]
		\centering
		\begin{subfigure}[b]{0.48\textwidth}
			\includegraphics[width=\textwidth]{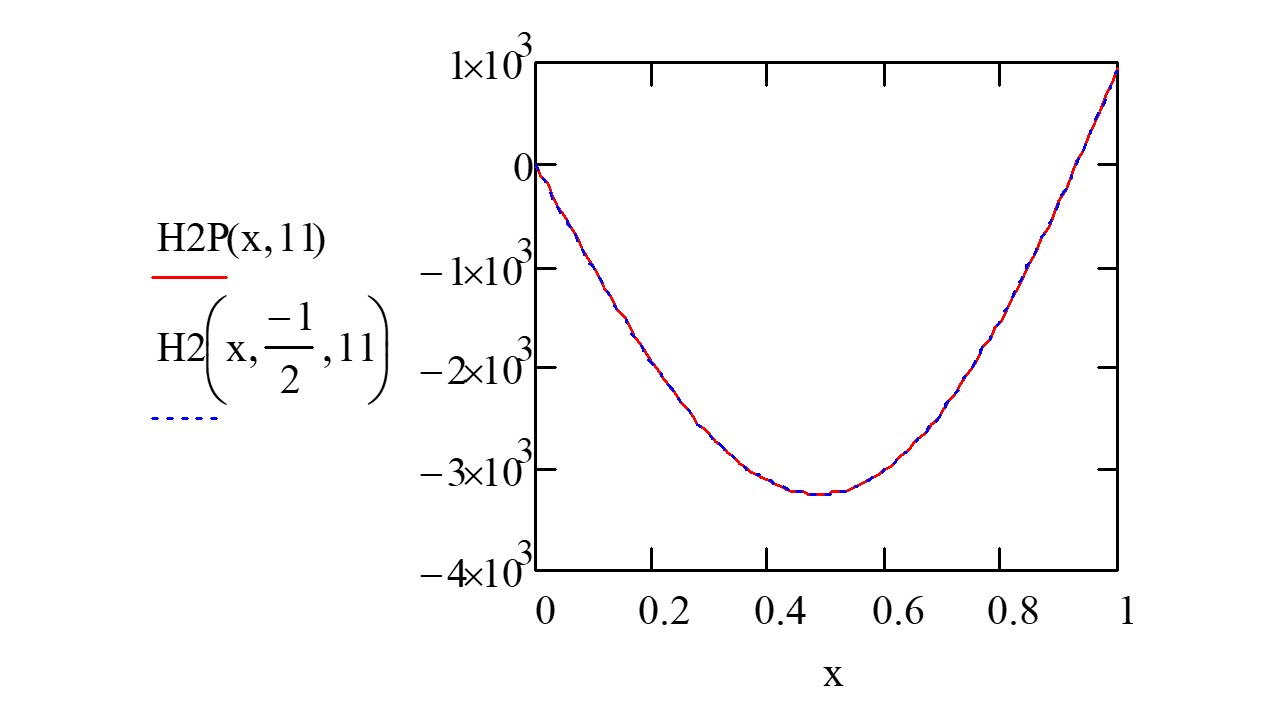}
			\caption*{\centering{\emph{Figure 4(a): Red $H2P(x,11)= {_{[3|2]}He_{11}(x)}$, Dot $H2\left(x,\dfrac{-1}{2},11\right)=He_{11}(x)$}}}
		\end{subfigure}
		\hfill
		\begin{subfigure}[b]{0.48\textwidth}
			\includegraphics[width=\textwidth]{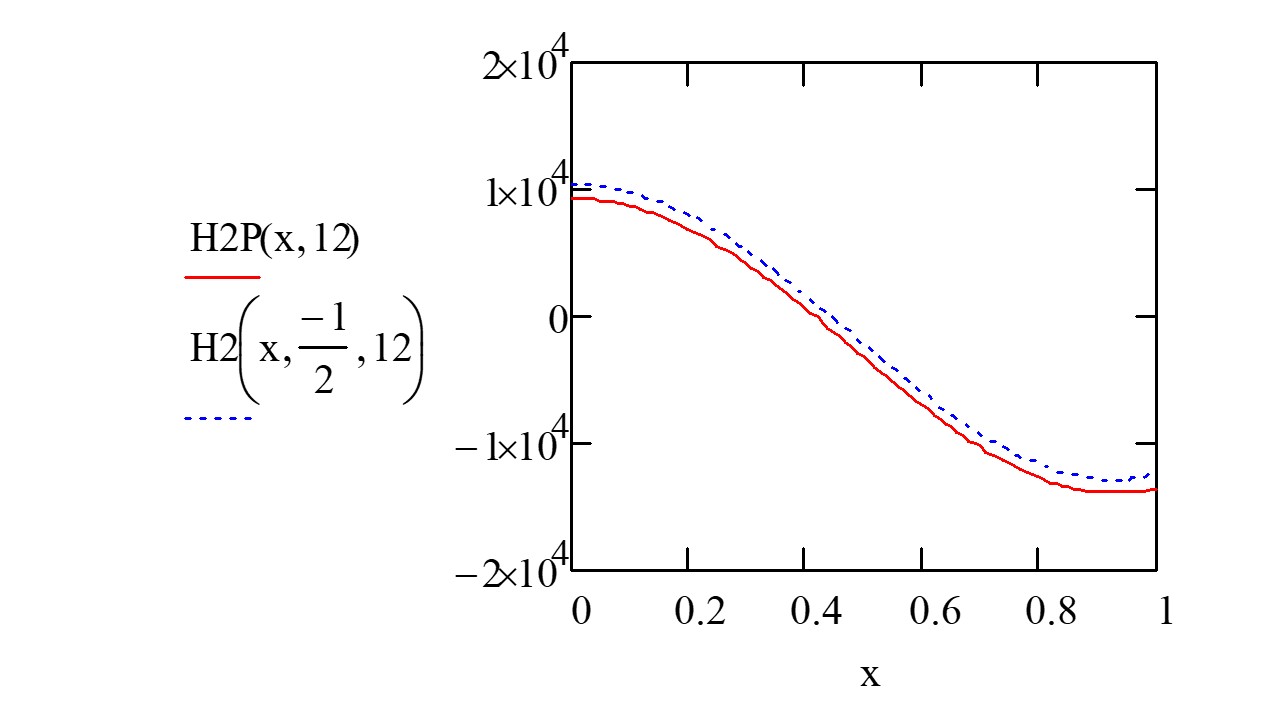}
			\caption*{\centering{\emph{figure 4(b): Red $H2P(x,12)= {_{[3|2]}He_{12}(x)}$, Dot $H2\left(x,\dfrac{-1}{2},12\right)=He_{12}(x)$}}}
		\end{subfigure}
		\centering
		\makebox[0.9\textwidth][c]{\em \textbf{Figure 3: Hermite polynomials $He_{n}(x)$ (dot) and their}}\\
		\makebox[0.9\textwidth][c]{\em \textbf{Pad\'e Approximation (red) for $n=11$ and $12$.}}
	\end{figure}
	
\noindent	The results of this section indicate that the method applied to the Appell polynomials has certain merits in revealing several interesting connections between them and other families of special polynomials.\\
	
The next section explores monomiality conditions and presents additional illustrative examples. Building on this foundation, the formalism is developed further to establish Padé approximants using the umbral notation.
	\section{Quasi monomials, Euler and Bernoulli polynomials and umbral formalism}
	This article has been devoted to the study of the consequences in terms of new families of polynomials ``generated'' by the Padè approximation of the amplitude of Appell polynomials. It is worth underscoring that the procedure, preserves the nature of the approximated polynomials, which are still recognized as belonging to the same polynomial set.\\
	
	The Appell polynomials are quasi monomials \cite{DG} an attribute, summarized through the following identities
	\begin{subequations}\label{paeq16}
		\begin{align}
			\hat{P} &= \partial_{x}, \label{paeq16a}\\
			\hat{M} &= x+ \dfrac{A^{\prime}(\partial_{x})}{A(\partial_{x})}. \label{paeq16b}
		\end{align}
	\end{subequations}
	such that
	\begin{subequations}\label{paeq17}
		\begin{align}
			\hat{P}a_{n}(x) &= na_{n-1}(x), \label{paeq17a}\\
			\hat{M}a_{n}(x) &=a_{n+1}(x). \label{paeq17b}
		\end{align}
	\end{subequations}
	A useful property can be deduced by using above operators, that is 
	\begin{equation}\label{paeq18}
		\hat{M}\hat{P}a_{n}(x)=na_{n}(x).
	\end{equation}
	The operators  $\hat{P}$ and $\hat{M}$ are known as derivative and multiplicative operators. Equation \eqref{paeq17} yields the recurrence properties and equation \eqref{paeq18} can be viewed as the differential equation satisfied by $a_n (x)$. In explicit terms, we find
	\begin{equation}\label{paeq19}
		\left[A(\partial_{x})x+A^{\prime}(\partial_{x})\right]\partial_{x}a_{n}(x)=nA(\partial_{x})a_{n}(x).
	\end{equation}
	The previous identity allows, for example, the derivation of differential equation characterizing the truncated exponential polynomials $e_{n}(x,1)=e_{n}(x)$. This can be easily obtained by noting that they are defined by the amplitude $A(t)=\dfrac{1}{1-t}$. The straightforward application of identity \eqref{paeq19} and a few algebraic manipulations yields the following second order ordinary differential equation with non-constant coefficients \cite{AND}:
	\begin{equation}\label{paeq78}
		xe_{n}^{\prime\prime}(x)-(x+n)e_{n}^{\prime}(x)=-ne_{n}(x).
	\end{equation}
Moreover, regarding the PA ${_{[1|1]}H_{n}^{(1)}(x,-1)}$ of Hermite polynomials $H_{n}^{(1)}(x,-1)$, the multiplicative and derivative operators are obtained by proving the following result:
	\begin{thm}
		For the second order Pad\'e approximant ${_{[1|1]}H_{n}^{(1)}(x,-1)}$ of Hermite polynomials, the following multiplicative and derivative operators exist:
		\begin{subequations}\label{paeq52}
			\begin{align}
				\hat{P}&=\partial_{x},\label{paeq52a}\\ \hat{M}&=x-\dfrac{1}{1-\dfrac{1}{4}\partial_{x}^{2}}.\label{paeq52b}
			\end{align}
		\end{subequations}
		\begin{proof}
			From equation \eqref{paeq17a}, derivative operator for ${_{[1|1]}H_{n}^{(1)}(x,-1)}$ is trivially $\hat{P}=\partial_{x}$.\\
			In view of equation \eqref{paeq11}, for $A(\partial_{x})=e^{-\partial_{x}}$, we have
			\begin{equation}\label{paeq20}
				{_{[1|1]}A(\partial_{x})}={_{[1|1]}e^{-\partial_{x}}}=\dfrac{\left(1-\dfrac{1}{2} \partial_{x}\right)}{\left(1+\dfrac{1}{2} \partial_{x}\right)},
			\end{equation}
		which on differentiating equation \eqref{paeq20} with respect to $x$ becomes
			\begin{equation}\label{paeq75}
				{_{[1|1]}A^{\prime}(\partial_{x})}=\dfrac{-1}{\left(1+\dfrac{1}{2}\partial_{x}\right)^{2} }.
			\end{equation}
		Use of equations \eqref{paeq20} and \eqref{paeq75} in equation \eqref{paeq16b} yields desired multiplicative operator \eqref{paeq52b}.
		\end{proof}
	\end{thm}
\begin{rem}
	Making use of operators \eqref{paeq52a} and \eqref{paeq52b} in monomiality equation \eqref{paeq18} for $(z_{n}(x):={_{[1|1]}H_{n}^{(1)}(x,-1)})$, the following result is deduced:
\end{rem}
	\begin{cor}
		The second order Pad\'e approximant of Hermite polynomials $z_{n}(x):={_{[1|1]}H_{n}^{(1)}(x,-1)}$ satisfy the following third order differential equation:
	\begin{equation}\label{paeq54}
	xz_{n}^{\prime\prime\prime}(x)+\left(2-n\right)z_{n}^{\prime\prime}(x)+4(1-x)z_{n}^{\prime}(x)+4nz_{n}(x)=0.
	\end{equation}
	\end{cor}
	Further examples regarding the explicit differential equations for the Pad\'e approximated polynomials will be discussed in the final section.\\
	
 The Euler and Bernoulli polynomials \cite{AND} have played an enormous role in pure and applied mathematics \cite{Milne}. They are members of the Appell set; accordingly, the considerations we have developed are suited for their study as well.
	
	The Euler polynomials $E_{n}(x)$ are specified by amplitude
	\begin{equation}\label{paeq29}
		A(t)=\dfrac{2}{e^{t}+1},\quad |t|<\pi.
	\end{equation} 
	Instead of treating the full amplitude approximation, for the sake of simplicity, we consider those with the exponential in the denominator. The use of the Taylor expansion of the exponential function in the denominator yields
\begin{equation}\label{paeq59}
			{_{[0|2]}A(t)}=\dfrac{1}{\left(1+\dfrac{1}{2}t+\dfrac{1}{4}t^{2}\right)},
	\end{equation}
	which in view of definition \eqref{paeq57} gives
	\begin{equation}\label{paeq76}
		{_{[0|2]}\left[\dfrac{2}{e^{t}+1}\right]}=\sum_{r=0}^{\infty}U_{r}\left(\dfrac{1}{2},\dfrac{1}{4}\right)t^{r}.
	\end{equation}
	Replacing $t$ by $\partial_{x}$ in equation \eqref{paeq76} and operating the resultant on $x^{n}$, it follows that
\begin{equation*}
	{_{[0|2]}\left[\dfrac{2}{e^{\partial_{x}}+1}\right]}\left\{x^{n}\right\}=\sum_{r=0}^{\infty}U_{r}\left(\dfrac{1}{2},\dfrac{1}{4}\right)\partial_{x}^{r}x^{n},
\end{equation*}
which accordingly yields
	\begin{equation}\label{paeq30}
		{_{[0|2]}E_{n}(x)}=n!\sum_{r=0}^{n}\dfrac{U_{r}\left(\dfrac{1}{2},\dfrac{1}{4}\right)}{(n-r)!}x^{n-r}.
	\end{equation}
	Higher-order approximants can be derived using a similar approach. The expression of the third-order PA ${_{[2|1]}E_{n}(x)}$ of Euler polynomials in terms of the  two-variable truncated exponential polynomials $e_{n}(x,y)$ is obtained in the following result:
	\begin{thm}
	The	third-order Pad\'e approximant $_{[2|1]}E_{n}(x)$ of Euler polynomials admits the following explicit representation in terms of the two-variable truncated exponential polynomials $e_{n}(x,y)$:
		\begin{equation}\label{paeq31}
			_{[2|1]}E_{n}(x) = e_{n}\left(x, -\dfrac{1}{12}\right) - \dfrac{5n}{12}e_{n-1}\left(x, -\dfrac{1}{12}\right) - \dfrac{n(n-1)}{24}e_{n-2}^{(2)}\left(x, -\dfrac{1}{12}\right).
		\end{equation}
		\begin{proof}
			The third-order PA ${_{[2|1]}A(t)}$ of the amplitude function $A(t) = \dfrac{2}{e^{t} + 1}$ is given by:
			\begin{equation*}
				_{[2|1]}A(t)= \dfrac{\left(1 - \dfrac{5}{12}t - \dfrac{1}{24}t^{2}\right)}{\left(1 + \dfrac{1}{12}t\right)}.
			\end{equation*}
		Replacing $t$ by $\partial_{x}$ and	substituting in operational definition \eqref{paeq10}, we arrive at the third-order PA $_{[2|1]}E_{n}(x)$ of the Euler polynomials:
			\begin{equation}\label{paeq65}
				_{[2|1]}E_{n}(x) = \dfrac{\left(1 - \dfrac{5}{12}\partial_{x} - \dfrac{1}{24}\partial_{x}^{2}\right)}{\left(1 + \dfrac{1}{12}\partial_{x}\right)}\{x^{n}\}.
			\end{equation}
			Rewrite equation \eqref{paeq65} as:
			\begin{equation*}
				_{[2|1]}E_{n}(x) = \left(\left(1 - \dfrac{5}{12}\partial_{x} - \dfrac{1}{24}\partial_{x}^{2}\right)\sum_{r=0}^{\infty}\left(-\dfrac{1}{12}\partial_{x}\right)^{r}\right)\left\{x^{n}\right\}.
			\end{equation*}
			Simplification of the above expression leads to desired result \eqref{paeq31}.
		\end{proof}
	\end{thm}
Further,	the polynomials 
	\begin{equation}\label{paeq66}
		{_{[0|2]}E_{n}\left(x;a,b\right)}=n!\sum_{r=0}^{n}\dfrac{U_{r}\left(a,b\right)}{(n-r)!}x^{n-r}
	\end{equation}
	can be understood as an umbral image \cite{BDLS} of truncated exponential polynomials, as reported below:
	\begin{equation}\label{paeq77}
		{_{[0|2]}E_{n}\left(x;a,b\right)}=e_{n}\left(x, \hat{u}\right)\varphi_{0},
	\end{equation}
	where 
		\begin{equation*}
			\hat{u}^{r}\varphi_{0}=U_{r}\left(a,b\right).
		\end{equation*}
	We have used the notation $E_{n}\left(x;a,b\right)$ instead of $E_{n}\left(x\right)$ here for the sake of clarity.\\
	
	It is to be noted that the polynomials in equation \eqref{paeq30} are a particular case of the polynomials in \eqref{paeq66} for $a=1/2$ and $b=1/4$.\\
	
	Similarly, the umbral representation for the third order PA $_{[1|2]}E_{n}(x;a,b)$ of the Euler polynomials is established in the below result:
	\begin{thm}
			For the third order Pad\'e approximant $_{[1|2]}E_{n}(x;1/6,1/12)$ of the Euler polynomials, the following umbral representation holds true:
		\begin{equation}\label{paeq68}
		_{[1|2]}E_{n}\left(x;\dfrac{1}{6},\dfrac{1}{12}\right)=e_{n}\left(x, \hat{u}\right)\varphi_{0}-\dfrac{n}{3}e_{n-1}\left(x, \hat{u}\right)\varphi_{0},
		\end{equation}
		where
		\begin{equation*}
			\hat{u}^{r}\varphi_{0}=U_{r}\left(\dfrac{1}{6},\dfrac{1}{12}\right).
		\end{equation*}
		\begin{proof}
		Proceeding on the same lines of proof as in Theorem \thmnumber{5}, the third-order PA $_{[1|2]}E_{n}(x)$ of the Euler polynomials is expressed as:
		\begin{equation*}
			_{[1|2]}E_{n}(x)=\dfrac{\left(1 - \dfrac{1}{3}\partial_{x}\right)}{\left(1+ \dfrac{1}{6}\partial_{x}+ \dfrac{1}{12}\partial_{x}^{2}\right)}\{x^{n}\}.
		\end{equation*}
	Using equation \eqref{paeq77} for $a=1/6$ and $b=1/12$, above equation takes the form
		\begin{equation}\label{paeq69}
					_{[1|2]}E_{n}\left(x;\dfrac{1}{6},\dfrac{1}{12}\right)=\left(1 - \dfrac{1}{3}\partial_{x}\right)e_{n}\left(x, \hat{u}\right)\varphi_{0},
		\end{equation}
		where
			\begin{equation*}
				\hat{u}^{r}\varphi_{0}=U_{r}\left(\dfrac{1}{6},\dfrac{1}{12}\right).
			\end{equation*}
		Simplification of equation \eqref{paeq69} yields assertion \eqref{paeq68}.
		\end{proof}
	\end{thm}
		In the previous discussion we have mentioned the Pad\'e approximants applied to the Euler polynomials. However, we did not explicitly examine the approximants of the exponential function appearing in the relevant Appell amplitude. The approximation, as defined in the introduction, is affordable without the need for computationally intensive algebraic steps. The umbral notation ensures a significant simplification, see, for example \cite{SLicciardi}.\\
		
		We note indeed that, regarding the Bernoulli polynomials $B_{r}(x)$ \cite{AND}, the amplitude of the relevant generating function has an exponential umbral image \cite{SUM}, namely 
		\begin{equation}\label{paeq44}
			A(t)=	\dfrac{t}{e^{t}-1}=e^{\hat{B}t}\phi_{0},\quad |t|<2\pi,\quad \hat{B}^{r}\phi_{0}=B_{r},
		\end{equation}
		with $B_{r}$  being the Bernoulli numbers \cite{AND} and the corresponding Bernoulli polynomials are defined as:
		\begin{equation}\label{paeq45}
			B_{r}(x)=(\hat{B}+x)^{r}\phi_{0}=H_{n}^{(1)}(x,\hat{B})\phi_{0}=\sum_{r=0}^{n}\binom{n}{r}B_{n-r}x^{r}.
		\end{equation}
		The Pad\'e approximation and its explicit representation can be derived by directly operating on the umbral image. As an illustration, we establish the second-order PA $_{[1|1]}B_{n}(x)$ of Bernoulli polynomials in the following result:
		
		\begin{thm}
			The second-order Pad\'e approximant $_{[1|1]}B_{n}(x)$ of Bernoulli polynomials admits the following umbral explicit representation:
			\begin{equation}\label{paeq47}
			{_{[1|1]}B_{n}(x)} = e_{n}\left(x,\dfrac{\hat{B}}{2}\right)\phi_{0} + \dfrac{n}{2}\hat{B}\;e_{n-1}\left(x,\dfrac{\hat{B}}{2}\right)\phi_{0}.
			\end{equation}
			\begin{proof}
				Starting with expression \eqref{paeq39}, we substitute $-t$ with $\hat{B}x$ and obtain the second-order PA ${_{[1|1]}A(t)}$ of the amplitude $A(t) = e^{\hat{B}t}\phi_{0}$:
				\begin{equation}\label{paeq46}
				{_{[1|1]}A(t)} = {_{[1|1]}e^{\hat{B}t}\phi_{0}} = \dfrac{\left(1 + \dfrac{\hat{B}}{2}t\right)}{\left(1 - \dfrac{\hat{B}}{2}t\right)}\phi_{0}.
				\end{equation}
				By incorporating above PA in operational definition \eqref{paeq10}, the second-order PA ${_{[1|1]}B_{n}(x)}$ of Bernoulli polynomials is expressed as:
				\begin{equation}\label{paeq70}
					{_{[1|1]}B_{n}(x)} = \dfrac{1 + \dfrac{\hat{B}}{2}\partial_{x}}{1 - \dfrac{\hat{B}}{2}\partial_{x}}\{x^{n}\}\phi_{0},
				\end{equation}
				which yields the result claimed in \eqref{paeq47}.
			\end{proof}
		\end{thm}
		Higher-order approximations can be obtained along the same lines by defining the higher-order PA through their umbral images. This approach can also be applied to the umbral forms of Euler and Genocchi polynomials \cite{DMS}.\\
		
		The results, obtained above, can be further extended. For example, consider the second-order Appell polynomials with the amplitude defined as $A(yt^{2})$. In this case, the corresponding Bernoulli polynomials, in umbral form writes
		\begin{equation}\label{paeq48}
			B_{n}^{(2)}(x,y)=H_{n}^{(2)}(x,y\hat{B})\phi_{0}.
		\end{equation}
		It is evident that, according to equation \eqref{paeq48}, the Hermite polynomials are the umbral image of the Bernoulli polynomials of order $2$ (or even higher). The relevant PA are therefore obtained by the direct use of the approximations derived in the previous sections.\\
		
		The forthcoming section contains some concluding comments on the directions along which this type of research may develop. In particular, we touch on the extension of the method to the case of Bessel functions. The Gaussian umbral interpretation of the $0^{th}$-order Bessel functions allows the derivation of interesting approximations.
		\section{Concluding remarks}
		This article has covered different aspects of the theory of Appell polynomials, viewed within the framework of the Padé approximant method. This point of view has been merged with more formal considerations, involving monomiality and umbral techniques.\\
		
	In order to outline a strategy for the study of their properties, we consider the case of two-variable truncated exponential polynomials $e_{n}(x,y)$, whose generating function is given in \eqref{paeq64}.
	The properties of these polynomials are easily derived from the previous definition itself. The recurrence under derivative with respect to the $x$-variable is straightforward, regarding that with respect to $y$-variable requires just a few comments. Keeping the $y$-derivative with respect to both sides of equation \eqref{paeq64}, we find
	\begin{equation}\label{paeq33}
		\partial_{y}e_{n}(x,y)=\dfrac{\partial_{x}}{1-y\partial_{x}}e_{n}(x,y),
	\end{equation}
	which eventually can be cast in the form
	\begin{equation}\label{paeq34}
		\partial_{y}e_{n}(x,y)=\sum_{r=0}^{n-1}\dfrac{n!y^{r}}{(n-1-r)!} e_{n-1-r}(x,y),
	\end{equation}
	where
	\begin{equation*}
		e_{n}(x,0)=x^{n}.
	\end{equation*}
	More in general, regarding the two-variable Appell polynomials $a_{n} (x,y)$ we can redefine the amplitude as $A(yt)$, and it comes out that the recurrence with respect to the $y$ derivative reads
	\begin{equation}\label{paeq35}
		\partial_{y}a_{n}(x,y)=T(y)\partial_{x} a_{n}(x,y);\quad a_{n}(x,0)=x^{n},
	\end{equation}
	where
	\begin{equation}\label{paeq40}
		T(y)=\dfrac{A^{\prime}(y\partial_{x})}{A(y\partial_{x})},
	\end{equation}
	with solution
	\begin{equation}\label{paeq71}
		a_{n}(x,y)=e^{\left(\int T(y)dy\right)\partial_{x}}\{x^{n}\}.
	\end{equation}
It is an evolutive partial differential equation reducing, in the case of Hermite polynomials, to the heat 
equation. The use of equation \eqref{paeq18} eventually yields 
\begin{equation}\label{paeq41}
	(y\partial_{y}+x\partial_{x})a_{n}(x,y)=(x+y\hat{T})\partial_{x}a_{n} (x,y)=na_{n} (x,y).
\end{equation}
The previous identity has a clear interpretation in terms of Lie operators \cite{Gil}. \\

The second order two-variable truncated exponetial polynomials $	e_{n}^{(2)}(x,y)$ defined as 
\begin{equation}\label{paeq42}
	e_{n}^{(2)}(x,y)=\dfrac{1}{1-y\partial_{x}^{2}}\{x^{n}\},
\end{equation}
which eventually leads to the following ordinary and partial differential equation
\begin{equation}\label{paeq43}	yxZ^{\prime\prime\prime}_{n}-nyZ^{\prime\prime}_{n}-xZ^{\prime}_{n}+nZ_{n}=0,\quad Z_{n}=e_{n}^{(2)}(x,y),
\end{equation}
	\begin{equation}\label{paeq72}
	\partial_{y}e_{n}^{(2)}(x,y)=\dfrac{\partial_{x}^{2}}{1-y\partial_{x}^{2}}e_{n}^{(2)}(x,y), \quad e_{n}^{(2)}(x,0)=x^{n},
\end{equation}
respectively.\\

Before closing the article, we like to emphasize the utility of umbral methods for the extension of the PA to other families of special functions. Therefore, we remind that the umbral image of the $0$-th order Bessel function $J_{0}(x)$ is a Gaussian function \cite{SLicciardi}, namely
\begin{equation}\label{paeq49}
	J_{0}(x)=e^{-\hat{c}\dfrac{x^{2}}{4}}\zeta_{0},
\end{equation}
where
\begin{equation*}
	\hat{c}^{n}\zeta_{0}=\dfrac{1}{\Gamma(n+1)}.
\end{equation*}
The use of second order PA ${_{[0|2]}J_{0}(x)}$ of Bessel function yields
\begin{equation}\label{paeq50}
	{_{[0|2]}J_{0}(x)}=\dfrac{1}{1+\dfrac{\hat{c}}{4}x^{2}+\left(\dfrac{\hat{c}}{4}\right)^{2}\dfrac{x^{4}}{2}}\zeta_{0},
\end{equation}
which eventually writes
\begin{equation}\label{paeq51}
	{_{[0|2]}J_{0}(x)}=\sum_{n=0}^{\infty}\dfrac{1}{r!}\left(\dfrac{x}{2}\right)^{2r}U_{r}\left(1,\dfrac{1}{2}\right).
\end{equation}
The validity of the approximation is checked in Figure $5(a)$, where we have reported the comparison of $J_{0}(x)$ with PA ${_{[0|2]}J_{0}(x)}$ and ${_{[0|3]}J_{0}(x)}$. As expected, the higher order approximants improve the agreement. The use of a fourth order PA allows a better matching in the $x$ interval from $0$ to $8$ as it can be seen in Figure $5(b)$.\\

\begin{figure}[h!]
	\centering
	\begin{subfigure}[b]{0.48\textwidth}
		\includegraphics[width=\textwidth]{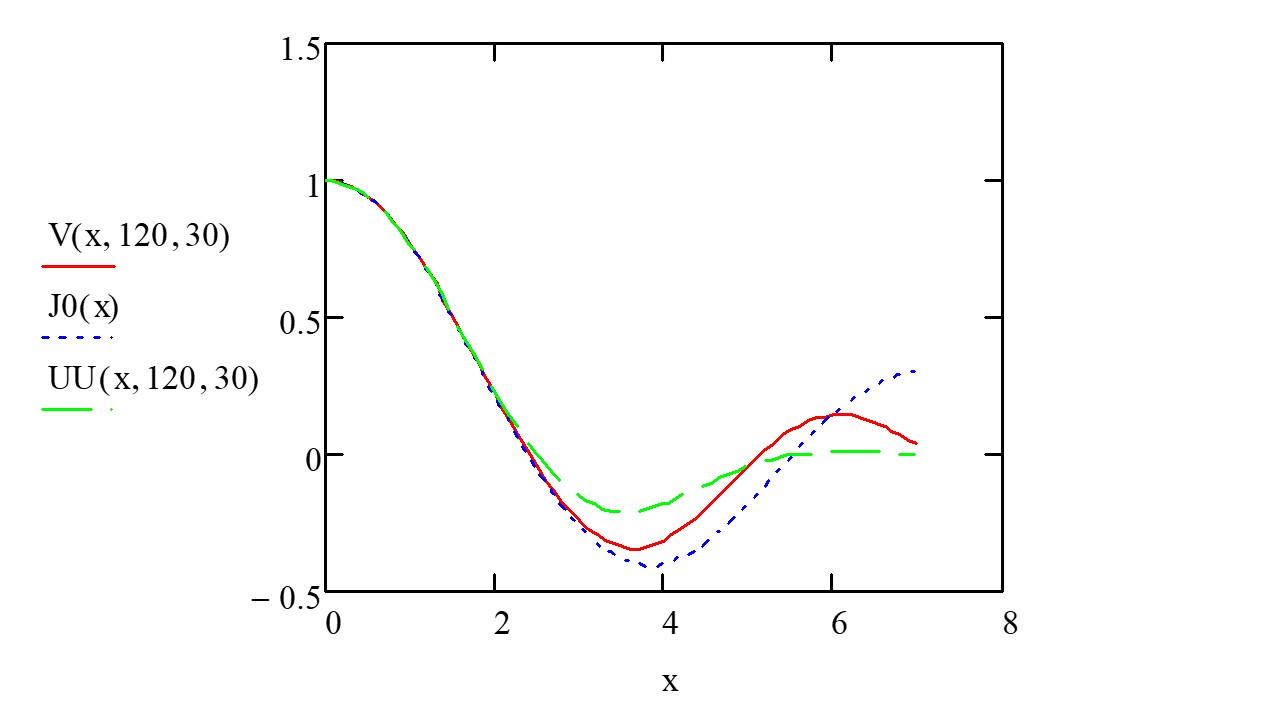}
		\caption*{\emph{Figure 5(a)}}
	\end{subfigure}
	\hfill
	\begin{subfigure}[b]{0.48\textwidth}
		\includegraphics[width=\textwidth]{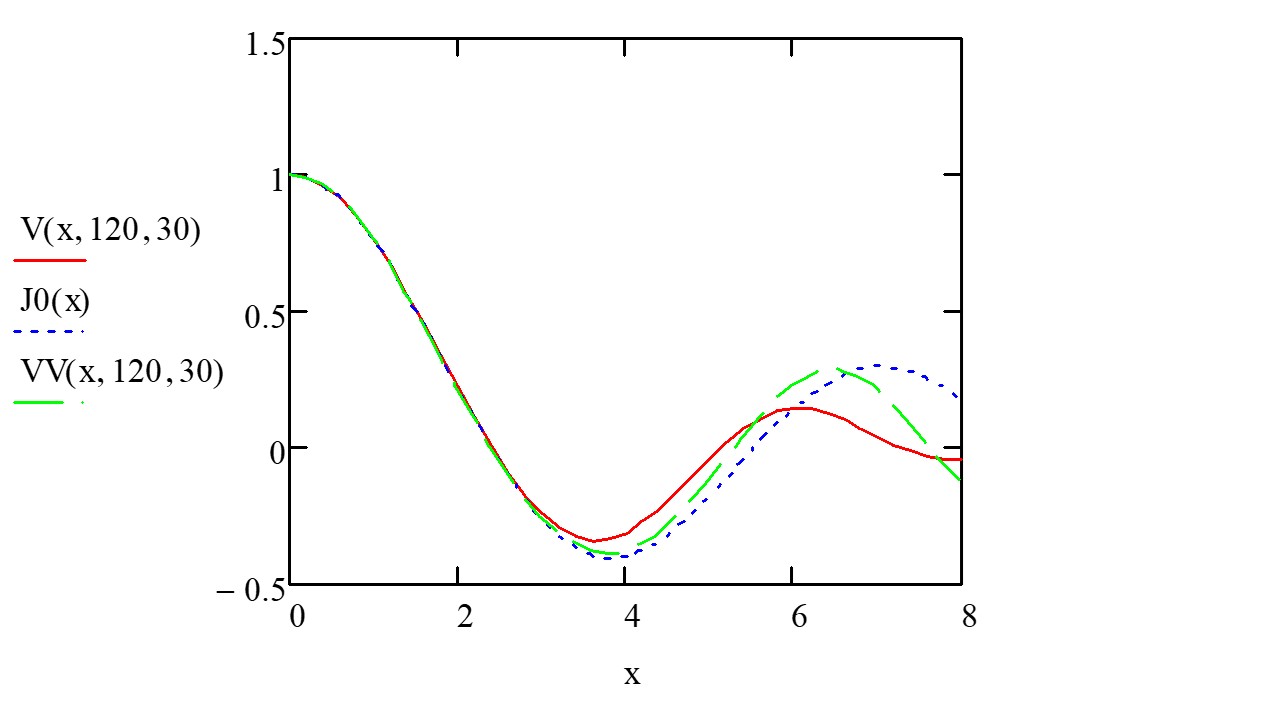}
		\caption*{\emph{Figure 5(b)}}
	\end{subfigure}
	{
		\centering
		\makebox[0.9\textwidth][c]{\em Figure 5(a) Comparison between the $0^{th}$ order cylindrical Bessel function (dot line)}\\		\makebox[0.9\textwidth][c]{ and Pad\'e approximants $[0|2]$ (solid line) and $[0|3]$ (dash line)}\\
	\makebox[0.9\textwidth][c]{\em Figure 5(b) Comparison between the $0^{th}$ order cylindrical Bessel function (dot line)}\\		\makebox[0.9\textwidth][c]{ and Pad\'e approximants $[0|2]$ (solid line) and $[0|4]$ (dash line)}}
\end{figure}
In this article, we have considered the theory of Appell polynomials. We have used the method of approximants merged with techniques of operational nature. The results we have obtained go beyond the mere approximation of polynomials. It seems to pave the way for a more general framing of the theory of Appell polynomials, which can be extended to the theory of Sheffer sequences \cite{Roman}.\\

\noindent {\bf Conflict of Interest:} The authors declare that they have no conflict of interest.

	
\end{document}